\documentclass[11pt]{article}
\usepackage{amssymb}
\usepackage{enumerate}
\usepackage{enumitem}
\usepackage{mathrsfs}
\usepackage{latexsym,amsmath,amssymb,amsfonts,epsfig,graphicx,cite,psfrag}
\usepackage{eepic,color,colordvi,amscd}

\newcommand{\qed}{\hfill $\Box $}
\newcommand{\pf}{\noindent {\bf Proof.} }
\newtheorem{theorem}{Theorem}
\newtheorem{lemma}[theorem]{Lemma}
\newtheorem{coro}[theorem]{Corollary}

\newtheorem{obser}[theorem]{Observation}

\begin{document}

\title{Spectral radius and rainbow matchings  of graphs
}

\author{
Mingyang Guo, Hongliang Lu\footnote{Corresponding email: luhongliang@mail.xjtu.edu.cn}, Xinxin Ma and Xiao Ma \\School of Mathematics and Statistics\\
Xi'an Jiaotong University\\
Xi'an, Shaanxi 710049, China\\
\smallskip\\
}
\date{}


\date{}

\maketitle

\begin{abstract}
Let $n,m$ be integers such that $1\leq m\leq (n-2)/2$ and let $[n]=\{1,\ldots,n\}$. Let $\mathcal{G}=\{G_1,\ldots,G_{m+1}\}$ be a family of graphs on the same vertex set $[n]$. In this paper, we prove that if for any $i\in [m+1]$, the spectral radius of $G_i$ is not less than $\max\{2m,\frac{1}{2}(m-1+\sqrt{(m-1)^2+4m(n-m)})\}$, then $\mathcal{G}$ admits a rainbow matching, i.e. a choice of disjoint edges $e_i\in G_i$, unless $G_1=G_2=\ldots=G_{m+1}$ and $G_1\in \{K_{2m+1}\cup (n-2m-1)K_1, K_m\vee (n-m)K_1\}$.
\end{abstract}

\section{Introduction}
For a set $S$, let ${S\choose 2}=\{\{x,y\} : \{x,y\}\subseteq S\}$.
A graph $G$ is a pair $G=(V,E)$, where $V:=V(G)$ is the vertex set and $E:=E(G)\subseteq {V(G)\choose 2}$ is the edge set. For $x,y\in V(G)$, we denote $\{x,y\}$ by $xy$ when there is no confusion.
A set of pairwise disjoint edges of $G$ is called a \emph{matching} in $G$. We use $\nu(G)$ to denote the maximum size of a matching in $G$.
Let $\mathcal{G}= \{G_1, ... ,G_t\}$ be a family of graphs on the same vertex set. A set of $t$ pairwise disjoint edges is called a \emph{rainbow matching} for $\mathcal{G}$ if each edge is from a different $G_i$. If such edges exist, then we also say that $\mathcal{G}$ \emph{admits} a rainbow matching.  Given a vertex subset $S\subseteq V(G)$, the subgraph induced by $S$ is denoted by $G[S]$. For a vertex $x\in V(G)$, let $N_G(x):=\{y\in V(G):xy\in E(G)\}$. We denote by $d_G(x)$ the number of edges containing $x$ in $G$.  The adjacency matrix of $G$ is $A(G)=(a_{ij})$, where $a_{xy}=1$ if $xy\in E(G)$ and $0$ otherwise. The eigenvalues of $G$ are the eigenvalues of its adjacency matrix $A(G)$. The largest eigenvalue of $A(G)$ is called the spectral radius of $G$ and denoted by $\rho(G)$. 


Let $G_1=(V_1,E_1)$ and $G_2=(V_2,E_2)$ be two graphs. The \emph{union} of  $G_1$ and $G_2$ denoted by $G_1\cup G_2$ is a graph  with vertex set $V_1\cup V_2$ and edge set $E_1\cup E_2$.
The \emph{join} of graphs $G_1$ and $G_2$ denoted by $G_1\vee G_2$ is a graph with vertex set $V_1\cup V_2$ and edge sets $E(G_1\cup G_2)$ together with all the edges joining $V_1$ and $V_2$.
We denote $G_1\simeq G_2$ if $G_1$ and $G_2$ are isomorphic. We say $G_1=G_2$ if $V(G_1)=V(G_2)$ and $E(G_1)=E(G_2)$.

For a positive integer $n$, let $[n]:=\{1,\ldots,n\}$ and let $[0]:=\emptyset$. For a graph $G$, we denote the complement of $G$ by $\overline{G}$. Let $m$ be an integer such that $1\leq m\leq (n-2)/2$. For each $1\leq i\leq m+1$, let $K^0_{i-1}$, $K^1_{2m-2i+3}$ and $K^2_{n-2m+i-2}$ denote complete graphs with vertex set $[i-1]$, $[2m-i+2]\setminus [i-1]$ and $[n]\setminus [2m-i+2]$ respectively and let $A^i_{n,m}:=K^0_{i-1}\vee (K^1_{2m-2i+3}\cup\overline{K^2_{n-2m+i-2}})$ .



Erd\H{o}s and Gallai \cite{EG59} determined the maximum $e(H)$ with $\nu(H)$ fixed. Akiyama and Frankl \cite{AF85} obtained the following result, which is a  generalization of Erd\H{o}s and Gallai's result \cite{EG59}.
\begin{theorem}[Akiyama and Frankl, \cite{AF85}]
Let $n, m$ be positive integers with $n\geq 2m+2$.
Let $\mathcal{G}= \{G_1, ... ,G_{m+1}\}$ be a family of graphs on the same vertex set $[n]$ such that $e(G_i)> \max\{{n\choose 2}-{n-m\choose 2}, {2m+1\choose 2}\}$ for each $i\in [m]$.
Then $\mathcal{G}$ admits a rainbow matching.
\end{theorem}

 Feng, Yu and Zhang \cite{FYZ17} gave a sufficient condition for a graph to have a matching in terms of spectral radius.
\begin{theorem}[Feng, Yu and Zhang, \cite{FYZ17}]\label{spectral}
	For any $n$-vertex graph $G$ with $\nu(G)\leq m$, we have
	\begin{enumerate}[itemsep=0pt,parsep=0pt]
		\item [$(1)$] if $n=2m$ or $2m+1$, then $\rho(G)\leq \rho(K_n)$ with equality if and only if $G=K_n$;
		\item [$(2)$] if $2m+2\leq n<3m+2$, then $\rho(G)\leq 2m$ with equality if and only if $G\simeq A^1_{n,m}$;
		\item [$(3)$] if $n=3m+2$, then $\rho(G)\leq 2m$ with equality if and only if $G\simeq A^{m+1}_{n,m}$ or $G\simeq A^1_{n,m}$;
		\item [$(4)$] if $n>3m+2$, then $\rho(G)\leq\frac{1}{2}(m-1+\sqrt{(m-1)^2+4m(n-m)})$ with equality if and only if $G\simeq A^{m+1}_{n,m}$.
	\end{enumerate}
\end{theorem}
Let $\mathcal{G}=\{G_1,\ldots,G_t\}$ be a family of not necessarily distinct graphs with common vertex set $V$. We say that a graph $H$ with vertex set $V$ is \emph{$\mathcal{G}$-rainbow} if there exists a bijection $\phi: E(H)\rightarrow[t]$ such that $e\in E(G_{\phi(e)})$ for each $e\in E(H)$. Joos and Kim \cite{JK20} considered the following question: Let $H$ be a graph with $t$ edges, $\mathbf{G}$ be a family of graphs and $\mathcal{G}=\{G_1,\ldots,G_t\}$ be collection of not necessarily distinct graphs on the same vertex set $V$ such that $G_i\in \mathbf{G}$ for all $1\leq i\leq t$. Which properties imposed on $\mathbf{G}$ yield a $\mathcal{G}$-rainbow graph isomorphic to $H$? In this paper,  we consider the property concerning spectral radius of graphs for the case $H$ is a matching of size $t$.

The following is our main result which is also a rainbow version of Theorem \ref{spectral}.
\begin{theorem}\label{main}
	Let $n, m$ be two positive integers such that $n\geq 2m+2$. Let $\mathcal{G}=\{G_1,G_2,\ldots,G_{m+1}\}$ be a family of  graphs on the same vertex set $[n]$. If
\begin{align}\label{main-eq}
 \rho(G_i)\geq\max\{2m,\frac{1}{2}(m-1+\sqrt{(m-1)^2+4m(n-m)})\},
\end{align}
 then
	\begin{enumerate}[itemsep=0pt,parsep=0pt]
		\item [$(1)$] for $2m+2\leq n<3m+2$, $\mathcal{G}$ admits a rainbow matching unless $G_1=\ldots=G_{m+1}\simeq A^1_{n,m}$;
		\item [$(2)$] for $n=3m+2$, $\mathcal{G}$ admits a rainbow matching unless $G_1=\ldots=G_{m+1}\simeq A^{m+1}_{n,m}$ or $G_1=\ldots=G_{m+1}\simeq A^1_{n,m}$;
		\item [$(3)$] for $n>3m+2$, $\mathcal{G}$ admits a rainbow matching unless $G_1=\ldots=G_{m+1}\simeq A^{m+1}_{n,m}$.
	\end{enumerate}
\end{theorem}

\section{Shifting}

In extremal set theory, one of the most important and widely-used tools is the technique of shifting, which allows us to limit our attention to sets with certain structure. Let $G$ be a graph with  vertex set $[n]$.  We define the \emph{$(x,y)$-shift} $S_{xy}$ by $S_{xy}(G)=\{S_{xy}(e):e\in E(G)\}$, where
$$S_{xy}(e)=\begin{cases}
	e\setminus \{y\}\cup \{x\}, & \text{if $y\in e$, $x\notin e$ and $e\setminus \{y\}\cup \{x\}\notin E(G)$;}\\
	e, &\text{otherwise.}
\end{cases}$$
If $S_{xy}(G)=G$ for every pair $(x,y)$ satisfying $x<y$ then $G$ is said to be \emph{shifted}.
\begin{obser}\label{obser}
If $G$ is shifted, then for any $\{x_1,x_2\},\{y_1,y_2\}\subseteq  [n]$ such that $x_i\leq y_i$ for $1\leq i\leq 2$, $\{y_1,y_2\}\in E(G)$ implies $\{x_1,x_2\}\in E(G)$.
\end{obser}
Iterating the $(x,y)$-shift for all pairs $(x,y)$ satisfying $x<y$ will eventually produce a shifted graph (see \cite{F87} and \cite{F95}). For a graph $G$, let $S(G)$ be a graph which is obtained from $G$ by the series of shifts and which is invariant under all possible shifts. Note that $S(G)$ is shifted. 

For graphs, the shifting operation is also known as the Kelmans operation (see \cite{BH12}). Csikv\'ari \cite{C09} proved that the spectral radius of a graph does not decrease under the Kelmans operation.

\begin{theorem}[Csikv\'ari, \cite{C09}]\label{Kelmans}
	Let $G$ be a graph on vertex set $[n]$. Let $u,v$ be two vertices of $G$. Then $\rho(S_{uv}(G))\geq \rho(G)$.
\end{theorem}

\begin{theorem}[Wu, Xiao and Hong, \cite{WXH}]\label{swich}
	Let $u,v$ be two vertices of the connected graph $G$. Suppose $v_1,v_2,\ldots,v_s\in N_G(v)\setminus (N_G(u)\cup \{u\})$ with $s\geq 1$ and $x=(x_1,x_2,\ldots,x_n)$ is the perron vector of $A(G)$, where $x_i$ corresponds to the vertex $v_i$ for $1\leq i\leq s$. Let $G^*$ be the graph obtained from $G$ by deleting the edges $\{v,v_i\}$ and adding the edges $\{u,v_i\}$ for $1\leq i\leq s$. If $x_u\geq x_v$, then $\rho(G^*)>\rho(G)$.
\end{theorem}

By applying Theorem \ref{swich}, we can prove a bit stronger result for connected graphs.
\begin{lemma}\label{shift-spectral}
Let $G$ be a connected graph on vertex set $[n]$. Let $u,v$ be two vertices of $G$. Then $\rho(S_{uv}(G))> \rho(G)$ unless $G\simeq S_{uv}(G)$.
\end{lemma}

\pf
Firstly, we consider $N_G(u)\backslash \{v\}\nsubseteq N_G(v)$ and $N_G(v)\backslash \{u\}\nsubseteq N_G(u)$. Let $\{x_1,\ldots,x_n\}$ be the perron vector of $A(G)$, where $x_u,x_v$ corresponds to the vertices $u,v$ respectively. For the case $x_u\geq x_v$, we have $\rho(S_{uv}(G))>\rho(G)$ by Theorem \ref{swich}. For the case $x_u< x_v$, we have $\rho(G)<\rho(S_{vu}(G))$ by Theorem \ref{swich}. Note that $S_{vu}(G)\simeq S_{uv}(G)$. Thus $\rho(S_{uv}(G))=\rho(S_{vu}(G))>\rho(G)$.

Secondly, we assume that $N_G(v)\backslash \{u\}\subseteq N_G(u)$ or $N_G(u)\backslash \{v\}\subseteq N_G(v)$. For the case $N_G(v)\backslash \{u\}\subseteq N_G(u)$, we have $S_{uv}(G)=G$. Thus $S_{uv}(G)\simeq G$. For the case $N_G(u)\backslash \{v\}\subseteq N_G(v)$, by the definition of shifting, one can see that $S_{uv}(G)\simeq G$.
\qed



\begin{lemma}[Huang, Loh and Sudakov, \cite{HLS}]\label{shift-matching}
	  Given a family of graphs $\mathcal{G}=\{G_1,\ldots,G_t\}$ on vertex set $[n]$. Let $x<y$ be two vertices in $[n]$. If $S_{xy}(\mathcal{G})$ admits a rainbow matching, then so does $\mathcal{G}$.
\end{lemma}	

\section{Proof of Theorem \ref{main}}

We prove Theorem \ref{main} by using the method in \cite{AF85} and some properties of spectral radius. The following lemma can be found in \cite{BH12}.
\begin{lemma}\label{Spec-incre}
If $u,v$ are two nonadjacent vertices of graph $G$, then $\rho(G+uv)>\rho(G)$.
\end{lemma}

\begin{lemma}\label{extremal-shift}
	Let $G$ be a graph on vertex set $[n]$. Let $x,y\in [n]$ such that $x<y$ and let $m$ be an integer such that $1\leq m\leq(n-2)/2$.
	\begin{enumerate}[itemsep=0pt,parsep=0pt]
		\item [$(1)$] If $\rho(G)=\rho(A^{m+1}_{n,m})$ and $S_{xy}(G)\simeq A^{m+1}_{n,m}$, then $G\simeq A^{m+1}_{n,m}$;
		\item [$(2)$] If $\rho(G)=\rho(A^1_{n,m})$ and $S_{xy}(G)\simeq A^1_{n,m}$, then $G\simeq A^1_{n,m}$.
	\end{enumerate}
\end{lemma}
	
\pf
Firstly, we prove (1). Since $S_{xy}(G)\simeq A^{m+1}_{n,m}$, there is a set $W$ of size $m$ such that $d_{S_{xy}(G)}(v)=n-1$ for every $v\in W$. 
If $\{x,y\}\subseteq W$ or $\{x,y\}\subseteq [n]\setminus W$, then $G=S_{xy}(G)\simeq A^{m+1}_{n,m}$. So we may assume that $|\{x,y\}\cap W|=1$.
Note that $d_{S_{xy}(G)}(x)\geq d_{S_{xy}(G)}(y)$. Hence we have $ \{x,y\}\cap W=\{x\}$ and so $d_{S_{xy}(G)}(x)=n-1$.  By the definition of $(x,y)$-shift, we can derive that $N_G(x)\cup N_G(y)=[n]$ and $xy\in E(G)$. It follows that $G$ is connected.
Since $\rho(G)=\rho(A^{m+1}_{n,m})=\rho(S_{xy}(G))$, we have $G\simeq S_{xy}(G)$ by Lemma \ref{shift-spectral}. Thus $G\simeq A^{m+1}_{n,m}$.

Secondly, we prove (2). Since $S_{xy}(G)\simeq A^1_{n,m}$, there is a set $U\subseteq [n]$ of size $2m+1$ such that $S_{xy}(G)[U]$ is a clique. If $|\{x,y\}\cap U|\neq 1$, then we have $G=S_{xy}(G)\simeq A^1_{n,m}$.
So we may assume that $|\{x,y\}\cap U|=1$. Then by the definition of $(x,y)$-shift, we have $x\in U$ and $y\notin U$. Let $G_1:=G[U\cup \{y\}]$ and $G_2:=S_{xy}(G)[U\cup \{y\}]$. Note that
\[
\rho(G_1)=\rho(G)=\rho(A^1_{n,m})=\rho(S_{xy}(G))=\rho(G_2).
\]
If $d_{G_1}(x)>0$ and $d_{G_1}(y)>0$, then $G_1$ is connected and $G_1, G_2$ are not isomorphic.  By Lemma \ref{shift-spectral}, we have $\rho(G_2)>\rho(G_1)$, a contradiction. Thus we have $N_{G_1}(x)=\emptyset$ or $N_{G_1}(y)=\emptyset$. It follows that $G\simeq S_{xy}(G)\simeq A^1_{n,m}$.  \qed

Recall that $S(G)$ is a graph which is obtained by the series of shifts and which is invariant under all possible shifts. The following result is an immediate consequence of Lemmas \ref{Kelmans} and \ref{extremal-shift} which will be used in the proof of Theorem \ref{main}.
\begin{coro}\label{extremal-shift2}
	Let $G$ be a graph on vertex set $[n]$ and let $m$ be an integer such that $1\leq m\leq (n-2)/2$.
	\begin{enumerate}[itemsep=0pt,parsep=0pt]
		\item [$(1)$] If $\rho(G)=\rho(A^{m+1}_{n,m})$ and $S(G)\simeq A^{m+1}_{n,m}$, then $G\simeq A^{m+1}_{n,m}$;
		\item [$(2)$] If $\rho(G)=\rho(A^1_{n,m})$ and $S(G)\simeq A^1_{n,m}$, then $G\simeq A^1_{n,m}$.
	\end{enumerate}
\end{coro}

\begin{lemma}\label{extremalcase}
Let $m,n$ be positive integers such that $n\geq 2m+2$ and $\mathcal{G}=\{G_1,G_2,\ldots,G_{m+1}\}$ be a family of graphs on the same vertex set $[n]$.  If one of the following statements holds, then $\mathcal{G}$ admits a rainbow matching.
	\begin{enumerate}[itemsep=0pt,parsep=0pt]
		\item [$(1)$] For each $1\leq i\leq m+1$, $G_i\simeq A^{m+1}_{n,m}$, and there exist $p,q\in [m+1]$ such that $W_p\neq W_q$, where $W_i=\{x\in [n] : d_{G_i}(x)=n-1\}$;
		\item [$(2)$] For each $1\leq i\leq m+1$, $G_i\simeq A^{1}_{n,m}$, and there exist $p,q\in [m+1]$ such that $U_p\neq U_q$, where $U_i=\{x\in [n] : d_{G_i}(x)=2m\}$.
	\end{enumerate}
\end{lemma}
\pf Firstly, suppose (1) holds. Let $\mathcal{W}=\{W_1,\ldots, W_{m+1}\}$.  Now we construct a bipartite graph $H$ with bipartition $(\mathcal{W},[n])$ such that $xW_i\in E(H)$ if and only if $x\in W_i$. Note that $W_p\neq W_q$ and $|W_i|=m$ for $1\leq i\leq m+1$.  Then we have $|\cup_{i\in S}W_i|\geq m$ for any $S\subseteq [m+1]$ and $|\cup_{i\in [m+1]}W_i|\geq m+1$.
So for any $X\subseteq \mathcal{W}$, we have $|N_{H}(X)|\geq |X|$. By Hall's Theorem, $H$ has a matching $M$ of size $m+1$. Without loss of generality, suppose that $M=\{x_1W_1,\ldots, x_{m+1}W_{m+1}\}$. By the definitions of $H$ and $W_i$, we have $d_{G_i}(x_i)=n-1$.

Write $X'=\{x_1,\ldots,x_{m+1}\}$. Since $n\geq 2m+2$, there are $m+1$ distinct vertices $y_1,\ldots,y_{m+1}\in V\setminus X'$. Recall that $d_{G_i}(x_i)=n-1$ for $1\leq i\leq m+1$. Thus $x_iy_i\in G_i$ for $1\leq i\leq m+1$ and $\{x_iy_i:1\leq i\leq m+1\}$ is a rainbow matching of $\mathcal{G}$.

Secondly, suppose (2) holds. Without loss of generality, we may assume that $p=1$ and $q=m+1$. Note that $U_1\neq U_{m+1}$. There exists a vertex $x_1\in U_{1}\setminus U_{m+1}$ such that $d_{G_1}(x_1)>0$. So we may find an edge $e_1=x_1y_1\in E(G_1)$. Now suppose we have found the matching $M_r=\{e_1,...,e_r\}$ such that $e_i=x_iy_i\in E(G_i-V(M_{i-1}))$ for $1\leq i\leq r$, where $M_0=\emptyset$. If $r=m+1$, then $M_{m+1}$ is a desired rainbow matching of $\mathcal{G}$. So we may assume that $r<m+1$. If $r=m$, since $U_1\neq U_{m+1}$ and $x_1\notin U_{m+1}$, we have
 \[
 |U_{m+1}\setminus V(M_m)|\geq (2m+1)-(2m-1)=2.
 \]
There exists an edge $x_{m+1}y_{m+1}\in E(G_{m+1}-V(M_m))$. Else if $r<m$, then
\[
 |U_{r+1}\setminus V(M_r)|\geq (2m+1)-(2m-2)= 3.
 \]
 So there exists an edge $x_{r+1}y_{r+1}\in E(G_{r+1}-V(M_r))$. Continuing this process for at most $m+1$ steps, we may obtain the desired rainbow matching.\qed

\medskip
\noindent \textbf{Proof of Theorem \ref{main}.}
Suppose that $\{G_1,\ldots,G_{m+1}\}$ dose not admit a rainbow matching. Then by Lemma \ref{shift-matching}, $\{S(G_1),\ldots,S(G_{m+1})\}$ dose not admit a rainbow matching. Let $H_i:=S(G_i)$ for $i\in [m+1]$. By the definition of shifting, Lemma \ref{Kelmans} and (\ref{main-eq}),  for $1\leq i\leq m+1$,
\begin{align}\label{mpr-eq1}
	\rho(H_i)\geq \rho(G_i)\geq \max\left\{2m,\frac{1}{2}(m-1+\sqrt{(m-1)^2+4m(n-m)})\right\}.
\end{align}

\medskip
\textbf{Claim 1.}
\begin{enumerate}[itemsep=0pt,parsep=0pt,label=$($\roman*$)$]
	\item  For $2m+2\leq n<3m+2$, $H_1= \ldots= H_{m+1}= A^1_{n,m}$;
	\item  For $n=3m+2$, $H_1= \ldots= H_{m+1}$ and $H_1\in \{A^1_{n,m},A^{m+1}_{n,m}\}$;
	\item  For $n>3m+2$, $H_1= \ldots= H_{m+1}= A^{m+1}_{n,m}$.
\end{enumerate}
\medskip

Let $e_i:=\{i,2m+3-i\}$ for $1\leq i\leq m+1$. We claim that $\{e_2,\ldots,e_m\}\subseteq E(H_i)$ for each $i\in [m+1]$. Otherwise, we assume that there exist $r\in [m]\backslash [1]$ and $t\in [m+1]$ such that $e_r=\{r,2m+3-r\}\notin E(H_t)$. Since $H_t$ is shifted, by Observation \ref{obser}, every edge $\{i,j\}\in E(H_t)$ satisfies $i<r$ or $j<2m+3-r$, where $i<j$. Thus $H_t$ is a subgraph of $A^r_{n,m}$.
Note that $\nu(A^r_{n,m})= m$ and $A^r_{n,m}$ is isomorphic to neither $A^1_{n,m}$ nor $A^{m+1}_{n,m}$ as $2\leq r\leq m$. By Theorem \ref{spectral} and Lemma \ref{Spec-incre}, we have
$$\rho(H_t)\leq \rho(A^r_{n,m})<\max\left\{2m,\frac{1}{2}(m-1+\sqrt{(m-1)^2+4m(n-m)})\right\},$$
a contradiction. Hence $\{e_2,\ldots,e_m\}\subseteq E(H_i)$ for each $i\in [m+1]$. Recall that $\{H_1,\ldots,H_{m+1}\}$ does not admit a rainbow matching. If there exist $p, q\in [m+1]$ such that $p\neq q$, $e_1\in H_p$ and $e_{m+1}\in H_q$, then $\{e_1,\ldots,e_{m+1}\}$ is a rainbow matching for $\{H_1,\ldots,H_{m+1}\}$, a contradiction. So we may distinguish three cases.

\medskip
\noindent\textbf{Case 1.} There exists $p\in [m+1]$ such that $e_1,e_{m+1}\in E(H_p)\setminus (\cup_{i\in [m+1]\setminus \{p\}}E(H_i))$.

Without loss of generality, we assume that $p=m+1$. Then $\{1,2m+2\},\{m+1,m+2\}\notin E(H_t)$ for each $t\in [m]$. Since $H_t$ is shifted, then by Observation \ref{obser},
  $\{1,2m+2\}\notin E(H_t)$ implies that $H_t$ is a subgraph of $A^1_{n,m}$. Since  $\{m+1,m+2\}\notin E(H_t)$,  one can see that  $H_t\neq A^1_{n,m}$. So by  Lemma \ref{Spec-incre}, we have
  $$\rho(H_t)<\max\left\{2m,\frac{1}{2}(m-1+\sqrt{(m-1)^2+4m(n-m)})\right\},$$ a contradiction.

\medskip
\noindent\textbf{Case 2.} $e_1\notin \cup_{i=1}^{m+1}E(H_i)$.

By the same argument as that of proof in Case 1, one can see that $H_i$ is a subgraph of $A^1_{n,m}$ for each $i\in [m+1]$. Note that $\rho(A^1_{n,m})=2m$. By Lemma \ref{Spec-incre} and (\ref{mpr-eq1}),
we have
\[
\rho(A^1_{n,m})\leq \max\{\frac{1}{2}(m-1+\sqrt{(m-1)^2+4m(n-m)}),2m\}\leq \rho(H_i)\leq\rho(A^1_{n,m}).
\]
Thus $\rho(A^1_{n,m})=\rho(H_i)=2m\geq \frac{1}{2}(m-1+\sqrt{(m-1)^2+4m(n-m)})$, which implies that $H_i=A^1_{n,m}$ for each $i\in[m+1]$ and $2m+2\leq n\leq 3m+2$.
%

\medskip
\noindent\textbf{Case 3.} $e_{m+1}\notin \cup_{i=1}^{m+1}E(H_i)$.

Since $H_t$ is shifted, by Observation \ref{obser}, $[n]\setminus[m]$ is an independent set in $H_i$. So $H_i$ is a subgraph of $A^{m+1}_{n,m}$ for each $i\in[m+1]$. By Lemma \ref{Spec-incre}, one can see that
\begin{align}\label{c3_eq1}
\rho(H_i)\leq \rho(A^{m+1}_{n,m})=\frac{1}{2}(m-1+\sqrt{(m-1)^2+4m(n-m)}).
\end{align}
Combining (\ref{mpr-eq1}) and (\ref{c3_eq1}), we have
\begin{align}\label{c3_eq2}
\rho(H_i)= \rho(A^{m+1}_{n,m})=\frac{1}{2}(m-1+\sqrt{(m-1)^2+4m(n-m)})
\end{align}
and
\[
\frac{1}{2}(m-1+\sqrt{(m-1)^2+4m(n-m)})\geq 2m,
\]
which implies $n\geq 3m+2$.
Recall that $H_i$ is a subgraph of $A^{m+1}_{n,m}$. By (\ref{c3_eq2}), we have $H_i= A^{m+1}_{n,m}$ for each $i\in[m+1]$. This completes the proof  of Claim 1.


By (\ref{mpr-eq1}) and Claim 1, we have
\begin{align*}
\rho(G_1)=\ldots=\rho(G_{m+1})=2m
\end{align*}
for $2m+2\leq n<3m+2$.
Thus by Corollary \ref{extremal-shift2}, we have $G_1\simeq\ldots\simeq G_{m+1}\simeq A^1_{n,m}$ for $2m+2\leq n<3m+2$.
If there exist $i,j\in [m+1]$ such that $G_{i}\neq G_j$, then by Lemma \ref{extremalcase}, $\{G_1,\ldots,G_{m+1}\}$ has a rainbow matching, a contradiction.
So we have $G_1=\cdots= G_{m+1}\simeq A^1_{m+1}$ for $2m+2\leq n<3m+2$. With similar discussion, we can derive that $G_1=\cdots= G_{m+1}\simeq A^1_{n,m}$ or $G_1=\cdots= G_{m+1}\simeq A^{m+1}_{n,m}$ for $n=3m+2$ and $G_1=\cdots= G_{m+1}\simeq A^{m+1}_{n,m}$ for $n>3m+2$. This completes the proof. \qed


\begin{thebibliography}{99}

	

\bibitem{AF85}
J. Akiyama and P. Frankl, On the size of graphs with complete-factors, \emph{J. Graph Theory}, \textbf{9} (1985), 197--201.

\bibitem{BH12}
A. Brouwer and W. Haemers,
\newblock Spectra of graphs,
\newblock {Universitext, Springer, New York, 2012}.

\bibitem{C09} P. Csikv\'ari, On a conjecture of V. Nikiforov, \emph{Discrete Math.}, \textbf{309} (2009), 4522--4526.

\bibitem{EG59}
P. Erd\H{o}s and T. Gallai, On the maximal paths and circuits of graphs,  \emph{Acta Math. Acad. Sci. Hungar.}, \textbf{10} (1959), 337--357.


\bibitem{FYZ17} L. Feng, G. Yu and X. Zhang, Spectral radius of graphs with given matching number, \emph{Linear Algebra Appl.}, \textbf{422} (2007), 133--138.

\bibitem{F87} P. Frankl, The shifting techniques in extremal set theory, In \emph{Surveys in Combinatorics}, Vol. 123 of \emph{London Mathematical Society Lecture Notes}, Cambridge University Press, (1987), pp. 81--110.

\bibitem{F95}
P.~Frankl,
\newblock Extremal set systems, in: Handbook of Combinatorics,
\newblock {Elsevier, Amsterdam, 1995}.



\bibitem{HLS}
H.~Huang, P.~Loh and B.~Sudakov,
\newblock The size of a hypergraph and its matching number,
\newblock {\em Combin. Probab. Comput.}, \textbf{21} (2012), 442--450.

\bibitem{JK20}
F. Joos and J. Kim,
\newblock On a rainbow version of Dirac's theorem,
\newblock {\em Bull. London Math. Soc.}, \textbf{52} (2020), 498--504.

\bibitem{WXH}
B. Wu, E. Xiao and Y. Hong,
\newblock The spectral radius of trees on $k$ pendant vertices,
\newblock {\em Linear Algebra Appl.}, \textbf{395} (2005), 343--349.

\end{thebibliography}
\end{document}